\documentclass[11pt,english,12p]{amsart}
\usepackage{amsmath}
\usepackage{amssymb}
\usepackage{tikz}
\usepackage[autostyle]{csquotes}
\usepackage{mathrsfs}
\usepackage{float}
\usepackage{tikz-cd}

\usepackage{pst-node}
\usetikzlibrary{patterns}
\usetikzlibrary{arrows,positioning,shapes,fit,calc}

 \usepackage[all,cmtip]{xy}

\usetikzlibrary{arrows,positioning,shapes,fit,calc}

\newtheorem{thm}{Theorem}[section]

\newtheorem{prop}[thm]{Proposition}

\newtheorem{defn}[thm]{Definition}
\newtheorem{exmp}[thm]{Example}

\newcommand{\R}{{\mathbb R}}

\newcommand{\C}{{\mathbb C}}
\newcommand{\Z}{{\mathbb Z}}

\selectfont

\title{A Minimal left ideal description of Geometric structures in dimensions $6$, $7$, and $8$}

\author[R. Su\'arez]{Ricardo Su\'arez}
\address{Departament of Mathematics, UNITO, Torino, Italy}
\email{ricardo.suarez532@csuci.edu}

\begin{document}
\maketitle

\normalsize

\begin{abstract} 
In this paper we relate minimal left ideals on Clifford algebras with special geometric structures in dimensions $6,7,$ and $8$.
\end{abstract}
\maketitle
\section{Introduction}
Spinorial descriptions of special geometric structures in dimensions $6,7,$ and $8$ have been recently explored in Riemannian geometry. The  spaces of spinors that are often used are $\Delta=\R^{m}$, where $m=8,8,16$, for $SU(3)$, $G_2$, and $Spin(7)$ structures respectively. The goal of this paper is to canonically identify special geometric structures in dimensions $6,7,$ and $8$  with the left Clifford modules generated by a primitive idempotents, i.e. a minimal left ideals for the appropriate dimensions. For instance, in dimension $6$, the $SU(3)$ structure $\psi_+,\psi_{-},\omega_{0}$ defines the minimal left ideal $\R_{0,6}\cdot \frac{1}{32}(\star q^{*}(\psi_{+}\wedge \psi_{-})+4q^{*}(\psi_{+})+4\star q^{*}(\omega_0))$, where $\frac{1}{32}(\star q^{*}(\psi_{+}\wedge \psi_{-})+4q^{*}(\psi_{+})+4\star q^{*}(\omega_0))$ is the unit element in the induced minimal left ideal.  In dimension $7$, the fundamental $3$-form $\phi$ that defines the $G_2$ structure induces the minimal left ideal 
$\R_{0,7}\cdot \frac{1}{112}(\star q^{*}(\phi\wedge \star \phi)+7 q^{*}(\phi)-7q^{*}(\star\phi)-q^{*}(\phi\wedge \star \phi))$, where $\frac{1}{112}(\star q^{*}(\phi\wedge \star \phi)+7 q^{*}(\phi)-7q^{*}(\star\phi)-q^{*}(\phi\wedge \star \phi))$ defines the  unit element in this minimal left ideal. Lastly, in dimension $8$, the $Spin(7)$ structure induced by the $4$-form $\Omega_{0}$ defines the minimal left ideal $\R_{0,8}\cdot \frac{1}{128}(\star q(\Omega_0\wedge \Omega_{0})-8q(\Omega_0)+q(\Omega_0\wedge\Omega_{0}))$.

\section{ Geometric structures in dimension $6,7$, and $8$}

  We begin with the definition of a $G$-structure on a differential manifold $M$.
\begin{defn}
Let $G$ be  a closed  Lie subgroup of $GL(n,\R)$. A \textbf{$G$-structure} on $M$ is a reduction of the structure group of the frame bundle $\mathcal{F}(M)$, $GL(n,\R)$, to $G$. That is, a $G$-structure is a  principal sub-bundle $Q$ of $\mathcal{F}(M)$ with fibre $G$. If $G$ is one of the groups that appear in the Berger classification theorem, then $G$ is called a \textbf{special geometric structure}. 
\end{defn}

Let $(M,g)$ be an oriented Riemannian spin manifold; that is, $(M,g)$ is a Riemannian manifold with vanishing first and second Stiefel Whitney classes. These manifolds admit at least two tensor fields, the Riemannian metric $g$ and the volume form $dV_{g}$, whose stabilizer group at any given point is the special orthogonal group $SO(n)$. As is well known for Riemannian manifolds the  holonomy group  $Hol(g)$ is strictly smaller than $SO(n)$  if and only if there  exist additional  nontrivial parallel tensor fields on $M$ with respect to  the Levi-Civita connection $\nabla^{g}$.  Thus  $G$-structures in general can be viewed in terms of the additional  tensor fields assigned to our oriented  Riemannian manifolds. In 1955 Marcel Berger classified the possible Lie subgroups of $ O(n)$ for  $n$-dimensional, simply connected manifolds  with an irreducible, non-symmetric Riemannian metric (see [14]).

\begin{thm}

Suppose that $M$ is a simply connected manifold of dimension $n$ and $g$ is an irreducible, non-symmetric Riemannian metric on $M$. Then exactly one of the following seven cases holds. 

\begin{enumerate}

\item $Hol(g)=SO(n)$.
\item $n=2m$ with $m\ge 2$, and $Hol(g)=U(m)$ in $SO(2m)$.
\item $n=2m$ with $m\ge 2$, and $Hol(g)=SU(m)$ in $SO(2m)$.
\item $n=4m$ with $m\ge 2$, and $Hol(g)=Sp(m)$ in $SO(4m)$.
\item $n=4m$ with $m\ge 2$, and $Hol(g)=Sp(m)Sp(1)$ in $SO(4m)$.
\item $n=7$  and $Hol(g)=G_2$ in $SO(7)$.
\item $n=8$ and $Hol(g)=Spin(7)$ in $SO(8)$.

\end{enumerate}
\end{thm}

 In dimension $6$, an $SU(3)$ structure can be defined by a $2$-form and  a stable $3$-form  which can  locally be expressed  as $\omega=e^{12}+e^{34}+e^{56}$ and $\psi_+=e^{135}-e^{146}-e^{236}-e^{245}$.

Moreover, on a $6$-dimensional oriented Riemannian manifold, we have two additional tensors, the  almost complex structure $J$ and a complex volume form $\Psi$ which can be expressed entirely in terms of $J$ and  $\psi_{+}$:  $\Psi=\psi_{+}+i\psi_{-}$, where $\psi_{-}=J\psi_{+}$. Moreover, the Riemannian metric can be expressed in terms of $\omega$ and  $J$, given by the tensor equation $g=\omega(\cdot,J\cdot)$.  In dimension $7$, a $G_2$ structure on a $7$-dimensional Riemannain manifold is defined by a stable $3$-form $\phi$ which can be locally (point-wise) defined as $\phi=e^{123}+e^{145}-e^{257}-e^{347}+e^{167}-e^{356}+e^{246}$ with respect to the local frame $e^{1},\ldots,e^{7}$ of $T_{p}^{*}M$. For $G_2$ manifolds, the Riemannian metric and volume form are expressed in terms of the  fundamental $3$-form via  the equation $g(X,Y)dV=\frac{1}{6}(\iota_X\phi)\wedge (\iota_Y\phi)\wedge \phi$  for any pair of vector fields $X,Y$. On dimension $8$ Riemannian manifolds, a $Spin(7)$ structure is defined as an admissible $4$-form that is locally expressed as: 
$$\Omega=e^{1234}+e^{1256}+e^{1278}+e^{1357}-e^{1368}-e^{1458}$$
$$-e^{1467}-e^{2358}-e^{2367}-e^{2457}+e^{2468}+e^{3456}+e^{3478}+e^{5678}.$$

\section{Real Clifford algebras and minimal left ideals generated by primitive idempotents }
Let $V$ be a finite dimensional $\R$ vector space with quadratic form $q$, and  $V^{\otimes}=\bigoplus_k V^{\otimes k}$ its tensor algebra. We define the \textbf{Clifford algebra} of  the pair $(V,q)$ as the quotient of the tensor algebra with the two-sided ideal $I_q=\langle v\otimes v+q(x)1_{V^{\otimes}}\rangle$; that is, $C_{q}(V)=V^{\otimes}/I_q$. The Clifford  algebra $C_{q}(V)$ carries a natural $\mathbb{Z}_2$ grading, where $C_{q}^{0}(V)$ denotes the elements of even degree, while $C_{q}^{1}(V)$ denotes the  elements of odd degree.  Choosing a basis for $(V,q)$, say $e_1,\ldots,e_n$,  we get the following canonical basis for the Clifford algebra generated by the following $2^n$ monomials:

\begin{equation}
\label{canonicalbasis}
\{1_V\} \cup \{e_{i_1} \cdots e_{i_k} : i_1 < i_2 < \cdots < i_k, k = 1, \ldots , n\}.
\end{equation}

As is well known, for any quadratic space over $\R$ we have an isomorphism with a quadratic space of the form $\R^{p,q}$, with a signature $(p,q)$, where $p$ is the number of positive definite generators and $q$ the number of negative definite generators. For  the quadratic spaces $\R^{p,q}$ and $\R^{n}$, we  denote the associated Clifford algebras by $\R_{p,q}$ and $\R_{0,n}$ respectively. $\R_{0,n}$ is canonically isomorphic to $\bigwedge \R^n$ as $\R$ vector spaces. This isomorphism is achieved via the assignment of the \textbf{quantization map} $q:\bigwedge \R^n\rightarrow \R_{0,n}$, where  $e_{i_1}\wedge \cdots \wedge e_{i_k}\mapsto e_{i_1} \cdots e_{i_k}$. 
The natural grading of both algebras is preserved via this assignment.
  The inverse of this isomorphism is what is called the \textbf{symbol map}, which we denote $\sigma: \R_{0,n}\rightarrow \bigwedge \R^n$, given by $\sigma(e_{i_1} \cdots e_{i_k})=e_{i_1}\wedge \cdots \wedge e_{i_k}$. These maps have natural extensions onto the spaces of exterior forms on $\R^n$. We denote these extensions by $q^{*}:\bigwedge^{*}(\R^n)^{*}\rightarrow \R_n$ where $q^{*}(e^{i_i}\wedge  \cdots \wedge e^{i_k})=e_{i_1} \cdots e_{i_k}$, and $\sigma^{*}:\R_n\rightarrow \bigwedge^{*} (\R^n)^{*}$ where $\sigma^{*}(e_{i_1} \cdots e_{i_k})=e^{i_1}\wedge \cdots \wedge e^{i_k}$, where $e^{1}, \ldots ,e^{n}$ is the dual basis to the canonical basis $e_1, \ldots ,e_n$; that is, $e^{i}(e_j)=\delta_{ij}$. It is well known that Clifford algebras are isomorphic to matrix algebras over the division algebras $\R,\C,$ and $\mathbb{H}$. This classification is given in the following theorem (see [7]).

\begin{thm}
\label{classtheorem}
The Clifford algebra $\R_{p,q}$, where $p+q=n$,  has the following minimal representations over $\R$, $\C$, and $\mathbb{H}$:

\begin{enumerate}
\item[i)] $\R_{p,q}\cong M_{2^{\frac{n}{2}}}(\R)$ if $q-p = 0,6 \mod 8$.

\item[ii)] $\R_{p,q}\cong M_{2^{\frac{n-1}{2}}}(\C)$ if $q-p = 1,5 \mod 8$.

\item[iii)] $\R_{p,q}\cong M_{2^{\frac{n-2}{2}}}(\mathbb{H})$ if $q-p = 2,4 \mod 8$.

\item[iv)] $\R_{p,q}\cong M_{2^{\frac{n-3}{2}}}(\mathbb{H})\oplus M_{2^{\frac{n-3}{2}}}(\mathbb{H})$ if $q-p = 3 \mod 8$.

\item[v)] $\R_{p,q}\cong M_{2^{\frac{n-1}{2}}}(\R)\oplus M_{2^{\frac{n-1}{2}}}(\R)$ if $q-p = 7 \mod 8$.

\end{enumerate}
\end{thm}
For any semisimple algebra $A$,  a minimal left ideal is of type $A\cdot e$, where $\cdot$ is multiplication in the algebra and $e\in A$ is a \textbf{primitive idempotent}. An element $e$ is primitive if it cannot be written as a sum of orthogonal  idempotents ($e=f+g$; $f^{2}=f$, $g^{2}=g$, $f\cdot g=g\cdot f=0$ see [17]).  An idempotent is called \textbf{minimal} if it is minimal  with respect to the partial ordering $f\le e$, which happens if and only if $ef=f=fe$ (see [17]). A \textbf{minimal left ideal} $I\subset A$  is a left ideal that does not contain any other nonzero left ideals. It is known that if $I$ is a minimal left ideal of our algebra $A$ then either $I^2=0$ or $I=A\cdot e$ for some idempotent $e\in A$ (see [17]). For our Clifford algebras, any minimal left ideal is of the form $\R_{p,q}\cdot f$ where $f$ is a primitive idempotent. $\R_{p,q}\cdot f$ is clearly a left $\R_{p,q}$ module, where module multiplication is given by $\R_{p,q}\times \R_{p,q}\cdot f\rightarrow \R_{p,q} \cdot f $, via $(\phi,\psi\cdot f)\mapsto (\phi\cdot \psi)\cdot f$, for all $\phi\in\R_{p,q}$ and $\psi\cdot f\in\R_{p,q}\cdot f.$

\begin{defn}
The minimal left ideals of the Clifford algebra  $\R_{p,q}$ are called spinor spaces, and the elements of the minimal left ideals are called algebraic spinors.

\end{defn}
As is consistent with spin geometry, these spinor spaces generate our spinor representations of the appropriate dimension. The following theorem gives us the construction and classification of minimal left ideals in $\R_{p,q}$ (see [17]).

\begin{thm}
A minimal left ideal of $\R_{p,q}$ is of type $\R_{p,q}\cdot f$ where $f=\dfrac{1+e_{t_1}}{2}\cdot \cdots \cdot \dfrac{1+e_{t_k}}{2}$ is a primitive idempotent in $\R_{p,q}$  and $e_{t_1}, \ldots ,e_{t_k}$ is a set of commuting elements of the canonical basis such that $e_{t_i}^2=1$ for all $i=1, \ldots ,k=q-r_{q-p}$; moreover, the generators form a multiplicative   group of order $2^{q-r_{q-p}}$. The numbers $r_i$ are called the Randon-Hurwitz numbers,  given by the recurrence formula $r_{q-p}$ subject to the conditions: $r_0=0,\ r_1=1,\ r_2=2,\ r_3=2,\ r_j=3$ where $4\le j\le 7$, $r_{i+8}=r_i+4$ for $i\ge 0$, $r_{-1}=-1$, and $r_{-i}=1-i+r_{i-2}$ for $i\ge 2$.

\end{thm}
The $k$ commuting elements generate $2^k$ different idempotents which yield the decomposition of the Clifford algebra given by: 
$$\R_{p,q}=\bigoplus_{all\ \pm\  combinations} \R_{p,q}\cdot \prod_{\alpha} \dfrac{1\pm e_{\alpha}}{2},$$ where each minimal left ideal 
$\R_{p,q}\cdot \prod_{\alpha} \dfrac{1\pm e_{\alpha}}{2}$ is of real dimension $2^{p+q-k}$. The algebra of endomorphisms  is isomorphic to the real matrix algebra of dimensions matching those of the above theorem; that is, $\R_{p,q}\cong End(\R_{p,q}\cdot f)\cong M_{2^{p+q-k}}(\R)$. Restricting our representations to the spin sub-groups  $\rho: Spin(p,q)\rightarrow Aut(\R_{p,q}\cdot f)$ gives us the usual spinor representations. In the next three sections, we properly recover geometric structures in dimensions $6,7,$ and $8$ in terms of these minimal left ideals.

\section{Recovering $SU(3)$ structures from Algebraic spinors in dimension $6$}

In this section we make the association between the local description of  the $SU(3)$ structure in dimension six and its associated minimal left ideal of $\R_{0,6}$.

\subsection{An $SU(3)$ structure recovered from a minimal left ideal  in dimension $6$ }
An $SU(3)$ structure in dimension $6$ is given by tensors $\omega_0, \psi_{\pm},J_0$ such that their stabilizer in $\R^6$ is the group $SU(3)$. For the Clifford algebra $\R_{0,6}$ we define the primitive idempotent $f=\left(\dfrac{1+e_{135}}{2}\right)\left(\dfrac{1-e_{146}}{2}\right)\left(\dfrac{1-e_{236}}{2}\right)=\frac{1}{8}(1+e_{135}-e_{146}-e_{236}-e_{245}-e_{3456}-e_{1234}-e_{1256})$,  where $\R_{0,6}\cdot f$ is the associated minimal left ideal. Now we  normalize the idempotent to get unit coefficients  which we denote  $W=8f=1+e_{135}-e_{146}-e_{236}-e_{245}-e_{3456}-e_{1234}-e_{1256}$ in $\R_{0,6}$. Utilizing the quantization and symbol maps, $q^{*}:\bigwedge (\R^6)^{*}\rightarrow \R_{0,6}$  and  $\sigma^{*}:\R_{0,6}\rightarrow \bigwedge (\R^6)^{*}$, as well as the grading projection map $\pi_{\alpha}: \R_{0,6}\rightarrow \R_{0,6}^{\alpha}$, where $\pi_{\alpha}(x)=\langle x \rangle_{\alpha}$ is the $\alpha$ graded component, we have $\langle W \rangle_0=1,\langle W \rangle_3=e_{135}-e_{146}-e_{236}-e_{245},$ and $\langle W \rangle_4=e_{3456}-e_{1234}-e_{1256}$. This gives us the graded decomposition  $W=\langle W \rangle_0+\langle W \rangle_3+\langle W \rangle_4$. Using the Clifford Hodge dual, we have the relationship $\star \langle W \rangle_3=e_{246}-e_{235}-e_{145}-e_{136}$, and $\star \langle W \rangle_4=-(e_{12}+e_{56}+e_{34})$. With these relations, we  recover the $SU(3)$ structure given by  $\sigma^{*}(\langle W \rangle_3)=\psi_+\in \bigwedge^{3} (\R^{6})^{*}$, $-\sigma^{*}(\star \langle W \rangle_3)=\psi_{-}\in \bigwedge^{3} (\R^{6})^{*}$, $-\sigma^{*}(\star \langle W \rangle_4)=\omega_{0}\in \bigwedge^{2} (\R^{6})^{*}$. 

\begin{prop}
Fix the  primitive idempotent $f=\frac{1}{8}(1+e_{135}-e_{146}-e_{236}-e_{245}-e_{3456}-e_{1234}-e_{1256})$  in $\R_{0,6}$ that defines the minimal left ideal $\R_{0,6}\cdot f$. Normalize the left ideal so that we have unit coefficients  in the sense 
$W=8f$. From this we recover the tensor fields that define  an   $SU(3)$ structure in dimension $6$. That is, $\sigma^{*}( \langle W \rangle_3)=\psi_+, -\sigma^{*}(\star \langle W \rangle_3)=\psi_{-}$, and  $-\sigma^{*}(\star \langle W \rangle_4)=\omega_{0}$.

\end{prop}
\subsection{Minimal left ideal  induced from an $SU(3)$ structure}
For the converse, we fix  an  $SU(3)$ structure in $\R^{6}$ defined by the local tensors  $\omega_0=e^{12}+e^{34}+e^{56}$ and $\psi_+=e^{135}-e^{146}-e^{236}-e^{245}$, and $\psi_{-}=J\cdot \psi_{+}$, where $J$ is the complex structure in $\R^{6}$. It is easy to see that 
$q^{*}(\psi_+)=e_{135}-e_{246}-e_{236}-e_{145}$, $q^{*}(\psi_{-})=e_{136}+e_{145}+e_{235}-e_{246}$, and $q^{*}(\omega_0)=e_{12}+e_{56}+e_{34}$. The volume form can be expressed in terms of $\psi_{+}$ and $\psi_{-}$ using the formula $\psi_{+}\wedge \psi_{-}=4e^{123456}$; hence the volume element of the Clifford algebra $\R_{0,6}$ can be expressed by the formula  $\frac{1}{4} q^{*}(\psi_{+}\wedge \psi_{-})=e_{123456}$. Now via the Clifford Hodge product we have $\frac{1}{4}\star q^{*}(\psi_{+}\wedge \psi_{-})=1$. Hence we have the formula $\frac{1}{32}(\star q^{*}(\psi_{+}\wedge \psi_{-})+4q^{*}(\psi_{+})+4\star q^{*}(\omega_0))= \frac{1}{8}(1+e_{135}-e_{146}-e_{236}-e_{245}-e_{3456}-e_{1234}-e_{1256})=\left(\dfrac{1+e_{135}}{2}\right)\left(\dfrac{1-e_{146}}{2}\right)\left(\dfrac{1-e_{236}}{2}\right)$, where $e_{135}, e_{146}, e_{236}$ are commuting positive definite elements in the canonical basis. Hence we have a primitive idempotent, establishing  the following:

\begin{prop}
The $SU(3)$ structure $\omega_0,\psi_{+},\psi_{-}$ in $\R^6$ determines the primitive idempotent  $f=\frac{1}{32}(\star q^{*}(\psi_{+}\wedge \psi_{-})+4q^{*}(\psi_{+})+4\star q^{*}(\omega_0))$ in the Clifford algebra $\R_{0,6}$, and the resulting minimal left ideal is $\R_{0,6}\cdot  \frac{1}{32}(\star q^{*}(\psi_{+}\wedge \psi_{-})+4q^{*}(\psi_{+})+4\star q^{*}(\omega_0))$.

\end{prop}
 Definitionally,  elements  of our minimal left ideal  $\sigma\in \R_{0,6}\cdot \frac{1}{32}(\star q^{*}(\psi_{+}\wedge \psi_{-})+4q^{*}(\psi_{+})+4\star q^{*}(\omega_0))$ are algebraic spinors generated by $SU(3)$ structures. Although it is enough to have this canonical identification, from a well defined $SU(3)$ structure $\psi_{+},\psi_{-},\omega_{0}$ we can always generate a primitive idempotent and hence a minimal left ideal in this canonical manner.

 Now the primitive idempotent, $f=\frac{1}{32}(\star q^{*}(\psi_{+}\wedge \psi_{-})+4q^{*}(\psi_{+})+4\star q^{*}(\omega_0))$, represents the identity class in the canonical projection  $\phi\mapsto \phi\cdot f$. Moreover, we have the $\R$  basis $f,e_2f,e_3f,e_5f,e_{23}f,e_{25}f,e_{35}f,e_{235}f$ in $\R_{0,6}\cdot f$, where each basis element is a unit spinor in the module.

\section{Algebraic spinors in dimension $7$}
In dimension $7$, a $G_2$ structure is a positive $3$-form whose local tensor in $\R^7$ is given by $\phi=e^{123}+e^{145}+e^{167}+e^{246}-e^{257}-e^{347}-e^{356}\in\bigwedge^{3}(\R^{7})^{*}$, and whose fundamental $4$-form with respect to the Hodge star operator is given by $\star \phi_0=e^{4567}+e^{2367}+e^{2345}+e^{1357}-e^{1346}-e^{1256}-e^{1247}\in \bigwedge^{4}(\R^7)^{*}$.

\subsection{ $G_2$ structure recovered from a minimal left ideal.}
We begin by fixing  the primitive idempotent $f=\frac{1}{16}(1+e_{123})(1+e_{145})(1-e_{257})(1+e_{167})$, which generates the  minimal left ideal $\R_{0,7}\cdot f$. We normalize the idempotent to obtain unit scalars via  $W=16f=1+e_{123}+e_{145}-e_{2345}-e_{257}-e_{1357}+e_{1247}-e_{347}+e_{167}-e_{2367}-e_{4567}-e_{1234567}+e_{1256}-e_{356}-e_{246}-e_{1346}$. Using the symbol and projection maps, we have $\langle W \rangle_3=e_{123}+e_{145}+e_{167}+e_{246}-e_{257}-e_{347}-e_{356}$, and $\langle W \rangle_4=e_{2367}-e_{4567}+e_{1346}+e_{1256}-e_{2345}-e_{1357}+e_{1247}$. In terms of the projection decomposition, we write $W=1+ \langle W \rangle_3+ \langle W \rangle_4-e_{1234567}$, where $e_{1234567}$ is the volume element in $\R_{0,7}$. Using the Clifford Hodge dual, we have  $\star \langle W \rangle_3= \langle W \rangle_4$ and $\star e_{1234567}=1$, and thus we have $W=\star e_{1234567}+ \langle W \rangle_3+\star  \langle W \rangle_3-e_{1234567}$. Hence, using the extended symbol map, we have $\sigma^{*}(\langle W \rangle_3)=\phi_{0}\in\bigwedge^{3} (\R^7)^{*}$, which is the desired $G_2$ structure in $\R^7$. The associated $4$-form that comes with a $G_2$ structure is given by $-\sigma^{*}(\langle W \rangle_4)=\star \phi_0$. From this we have the following. 

\begin{prop}
Fix the primitive idempotent $f=\frac{1}{16}(1+e_{123})(1+e_{145})(1-e_{257})(1+e_{167})$, with  $\R_{0,7}\cdot f$ the associated  minimal left ideal. We  can then obtain a $G_2$ structure in $\R^7$ generated from the normalized tensor $W$ via $\sigma^{*}(\langle W \rangle_3)=\phi_0$, and the associated $4$-form is then given by $-\sigma^{*}(\langle W \rangle_4)=\star \phi_0$. 
\end{prop}

\subsection{Algebraic spinors generated by $G_2$ structures in dimension $7$}

Fix a $G_2$ structure  $\phi=e^{123}+e^{145}+e^{167}+e^{246}-e^{257}-e^{347}-e^{356}\in\bigwedge^{3}(\R^{7})^{*}$. Using the quantization map, we have  $q^{*}(\phi)=e_{123}+e_{145}+e_{167}+e_{246}-e_{257}-e_{347}-e_{356}$ and $q^{*}(\star\phi)=-e_{2367}+e_{4567}-e_{1346}-e_{1256}+e_{2345}+e_{1357}-e_{1247}$. The volume form in $\R^7$ expressed in terms of the $G_2$ structure is given by the formula $\phi\wedge \star\phi=7 e^{1234567}$, and thus  we have on $\R_{0,7}$ the equation  $q^{*}(\phi\wedge \star \phi)=7e_{1234567}$. Using the Clifford Hodge star operator, we have $\star q^{*}(\phi\wedge \star\phi)=7$. Putting this all together, we get the following primitive idempotent element in $\R_{0,7}$ induced by the $G_2$ structure: 
$$f_{\phi}=\frac{1}{112}(\star q^{*}(\phi\wedge \star \phi)+7 q^{*}(\phi)-7q^{*}(\star\phi)-q^{*}(\phi\wedge \star \phi))=\frac{1}{16}(1+e_{123})(1+e_{145})(1-e_{257})(1+e_{167}).$$

\begin{prop}
Fix the local  $G_2$ structure $\phi=e^{123}+e^{145}+e^{167}+e^{246}-e^{257}-e^{347}-e^{356}$ in $\R^7$. The $3$-form $\phi$ then determines the primitive idempotent  $f_{\phi}=\frac{1}{112}(\star q^{*}(\phi\wedge \star \phi)+7 q^{*}(\phi)-7q^{*}(\star\phi)-q^{*}(\phi\wedge \star \phi))$   in the Clifford algebra $\R_{0,7}$, with the resulting minimal left ideal  $\R_7\cdot \frac{1}{112}(\star q^{*}(\phi\wedge \star \phi)+7 q^{*}(\phi)-7q^{*}(\star\phi)-q^{*}(\phi\wedge \star \phi))=:\R_7\cdot f_{\phi}$.

\end{prop}

Elements of the minimal left ideal $\R_{0,7}\cdot f_{\phi}$ are algebraic spinors generated by the  $G_2$ structure. Although it is enough to have this canonical identification, from a well defined $G_2$ structure $\phi$ we can always generate a primitive idempotent and hence a minimal left ideal in this canonical manner.  As vector spaces we have  $\R_{0,7}\cdot f_{\phi}\cong \R^8$, with the basis of equivalence classes given by $f_{\phi},e_1f_{\phi},e_2f_{\phi},e_3f_{\phi},e_{4}f_{\phi},e_{5}f_{\phi},e_{6}f_{\phi},e_{7}f_{\phi}$.

\section{Algebraic spinors in dimension $8$}
In dimension $8$, a $Spin(7)$ structure is defined by the model tensor   $\Omega_{0}=e^{1234}+e^{1256}+e^{1278}+e^{1357}-e^{1368}-e^{1458}-e^{1467}-e^{2358}-e^{2367}-e^{2457}+e^{2468}+e^{3456}+e^{3478}+e^{5678}\in\wedge^{4} (\R^{8})^{*}.$ Using the quantization map, we have the following:
$q^{*}(\Omega_0)=e_{1234}+e_{1256}+e_{1278}+e_{1357}-e_{1368}-e_{1458}-e_{1467}-e_{2358}-e_{2367}-e_{2457}+e_{2468}+e_{3456}+e_{3478}+e_{5678},$   $q^{*}(\Omega_0\wedge \star\Omega_0)=8e_{12345678}$, and $\star q^{*}(\Omega_0\wedge \star\Omega_0)=8$. Now the formula $$f_{\Omega}=\frac{1}{128}(\star q(\Omega_0\wedge \Omega)-8q(\Omega_0)+q(\Omega_0\wedge\Omega_{0}))$$  is a primitive idempotent in $\R_{0,8}$, as it factors out as $f_{\Omega}=\frac{1}{16}(1-e_{1234})(1-e_{1256})(1-e_{1278})(1-e_{1357})$. Thus $\R_{0,8}\cdot\frac{1}{128}(\star q(\Omega_0\wedge \Omega)-8q(\Omega_0)+q(\Omega_0\wedge\Omega_{0}))$ is the minimal left ideal induced by our $Spin(7)$ structure. Moreover, $\R_{0,8}\cdot \frac{1}{128}(\star q(\Omega_0\wedge \Omega)-8q(\Omega_0)+q(\Omega_0\wedge\Omega_{0}))\cong\Delta_{8}=\R^{16}$.  We summarize this with the following proposition.

\begin{prop}
The $Spin(7)$ structure $\Omega_0$ in $\R^8$ determines the primitive idempotent $f_{\Omega}=\frac{1}{128}(\star q(\Omega_0\wedge \Omega)-8q(\Omega_0)+q(\Omega_0\wedge\Omega_{0}))$  in the Clifford algebra $\R_{0,8}$, and the resulting minimal left ideal is  $\R_{0,8}\cdot \frac{1}{128}(\star q(\Omega_0\wedge \Omega)-8q(\Omega_0)+q(\Omega_0\wedge\Omega_{0}))$.

\end{prop}

Conversely, for $\R_{0,8}$ we have four commuting generators from our canonical basis $e_{1234}, e_{1256}, e_{1278}, e_{1357}$, in which the idempotent $f=\frac{1}{16}(1-e_{1234})(1-e_{1256})(1-e_{1278})(1-e_{1357})$  defines  the spinor space $\R_{0,8}\cdot f$ isomorphic to $\Delta_{8}=\R^{16}$. The normalized tensor,  $W=16f=1-(e_{1234}+e_{1256}+e_{1278}+e_{1357}-e_{1368}-e_{1458}-e_{1467}-e_{2358}-e_{2367}-e_{2457}+e_{2468}+e_{3456}+e_{3478}+e_{5678})+e_{12345678}$, decomposes as $W=\star e_{12345678}-\langle W \rangle_4+e_{12345678}$. Hence we define $\sigma^{*}(\langle W \rangle_4)=\Omega_0$, and $\sigma^{*}(\star\langle W \rangle_4)=\star\Omega_0=\Omega_0$. This results in  the following proposition:

\begin{prop}
 Fix the primitive idempotent $f=\frac{1}{16}(1-e_{1234})(1-e_{1256})(1-e_{1278})(1-e_{1357})$ such that $\R_{0,8}\cdot f$ is the minimal left ideal.  We can associate a $Spin(7)$ structure in $\R^8$ generated from the normalized tensor $W$ via $\sigma^{*}(\langle W \rangle_4)=\sigma^{*}(\star\langle W \rangle_4)=\Omega_0$. 
\end{prop}

\section{Primitive idempotents in dimension $6$ to $G_2$ structures in dimension $7$}

We conclude this paper by relating the constructions in dimension $6$ and dimension $7$ by viewing  $\R^{7}=\R^{6}\oplus \R$, where the orthogonal dimension is given by the basis vector $e_7$. We define a  generic $SU(3)$  structure given in the $\R^6$ component  by  $\omega_0=e^{12}+e^{34}+e^{56}$, $\psi_{+}=e^{135}-e^{146}-e^{236}-e^{245}$, and $\psi_{-}=e^{136}+e^{145}+e^{235}-e^{246}$. For the associated Clifford algebra $\R_{0,6}$ we have the associated primitive idempotent $f=\left(\dfrac{1+e_{135}}{2}\right)\left(\dfrac{1-e_{146}}{2}\right)\left(\dfrac{1-e_{236}}{2}\right)=\frac{1}{8}(1+e_{135}-e_{146}-e_{236}-e_{245}-e_{3456}-e_{1234}-e_{1256})$,   with $\R_{0,6}\cdot f$  being the associated minimal left ideal. As we saw above from the normalized idempotent $W$ in dimension $6$,  we recover the $SU(3)$ structure via $\sigma^{*}(\langle W \rangle_3)=\psi_{+},\sigma^{*}(\star\langle W \rangle_3)=\psi_{-}$, and $\sigma^{*}(\star \langle W \rangle_4)=\omega_0$.  Now from the induced $SU(3)$ structure in dimension $6$, we can recover a $G_2$ structure in dimension $7$  via  $ \phi_0=\sigma^{*}(\star\langle W \rangle_4)\wedge e^7+\sigma^{*}(\langle W \rangle_3)$.  Conversely,  starting with the $SU(3)$ we have the primitive idempotent   $\frac{1}{32}(\star q^{*}(\psi_{+}\wedge \psi_{-})+4q^{*}(\psi_{+})+4\star q^{*}(\omega_0))$, which gives us the minimal left ideal $\R_{0,6}\cdot \frac{1}{32}(\star q^{*}(\psi_{+}\wedge \psi_{-})+4q^{*}(\psi_{+})+4\star q^{*}(\omega_0))$. Now the normalized idempotent with unit scalars $W$ in this formulation is given by $W=\star q^{*}(\psi_{+}\wedge \psi_{-})+4q^{*}(\psi_{+})+4\star q^{*}(\omega_0)$,  and hence we define the primitive idempotent in dimension $7$ as follows:
$f_{\phi}= \frac{1}{112}(\star q^{*}(\sigma^{*}(\star\langle W \rangle_4)\wedge e^7+\sigma^{*}(\langle W \rangle_3)\wedge \star \sigma^{*}(\star\langle W \rangle_4)\wedge e^7+\sigma^{*}(\langle W \rangle_3))+7 q^{*}(\sigma^{*}(\star\langle W \rangle_4)\wedge e^7+\sigma^{*}(\langle W \rangle_3))-7q^{*}(\star\sigma^{*}(\star\langle W \rangle_4)\wedge e^7+\sigma^{*}(\langle W \rangle_3))-q^{*}(\sigma^{*}(\star\langle W \rangle_4)\wedge e^7+\sigma^{*}(\langle W \rangle_3)\wedge \star \sigma^{*}(\star\langle W \rangle_4)\wedge e^7+\sigma^{*}(\langle W \rangle_3)
)).$  Thus $\R_{0,7}\cdot f_{\phi}$ is the spinor module induced from the $SU(3)$ structures recovered from the primitive idempotents in dimension $6$.

\section{Future research} 
With the correspondences established above, we can construct algebraic spinor bundles from the induced spinor spaces and provide classification equations for $G_2$, $SU(3)$, and $Spin(7)$ manifolds in terms of algebraic spinors.

\section{Acknowledgments }
I would like to thank my advisor, Dr. Ivona Grzegorczyk, and my supervisor at UNITO, Dr. Anna Fino, as well as  UNITO for financing my research and providing me guidance in this topic.

\end{document}